\documentclass{article}
\usepackage{amssymb}

\newtheorem{theorem}{Theorem}

\newtheorem{lemma}[theorem]{Lemma}

\input{tcilatex}
\begin{document}

\title{Exact inversion of Funk-Radon transforms with non-algebraic geometries%
}
\author{Victor Palamodov}
\maketitle

\textbf{Abstract.} Any even function defined on 2-sphere is reconstructed
from its integrals over big circles by means of the classical Funk formula.
For the non-geodesic Funk transform on the sphere of arbitrary dimension,
there is the explicit inversion formula\ similar to that for the geodesic
transform. A function defined on the sphere of radius one is integrated over
traces of hyperplanes tangent to a sphere contained in the unit ball. This
reconstruction is generalized in the paper for Riemannian hypersurfaces in
an affine space.

\textbf{MSC (2010) }Primary 53C65; Secondary 44A12

\section{Introduction}

A Riemannian manifold $X$ is embedded as a hypersurface in an affine space.
A function defined on $X$ is integrated over intersections of $X$ with
hyperplanes tangent to\ an ellipsoid $\Sigma \ $(called \textit{cam}). We
prove the reconstruction formula that looks like the inversion formula for
the non-geodesic Funk transform on the sphere stated in \cite{Pa16}. The
only condition on $X$ and $\Sigma $ is that no three points $x,y,\sigma $
are collinear.

\section{Auxiliary}

Let $X$ and $\Sigma $ be manifolds of dimension $n>1$ with volume forms $%
\mathrm{d}X$ and $\mathrm{d}\Sigma $ and $\Phi $ be a real smooth function
defined on $X\times \Sigma $ such that $\mathrm{d}\Phi \left( x,\sigma
\right) \neq 0$ as $\Phi \left( x,\sigma \right) =0.$ The Funk-Radon
transform $\mathrm{M}_{\Phi }$ generated by this function is defined by 
\[
\mathrm{M}_{\Phi }f\left( \sigma \right) =\lim_{\varepsilon \rightarrow 0}%
\frac{1}{2\varepsilon }\int_{\left\vert \Phi \right\vert \leq \varepsilon }f%
\mathrm{d}X=\int_{Z\left( x\right) }f\left( x\right) \frac{\mathrm{d}X}{%
\mathrm{d}_{x}\Phi },\ \sigma \in \Sigma .
\]%
Suppose that (\textbf{I}):\ the\textit{\ }map$\ \mathrm{D}:Z\times \mathbb{R}%
_{+}\rightarrow T^{\ast }\left( X\right) \backslash 0$\textit{\ }is a
diffeomorphism, where $Z=\left\{ \Phi \left( x,\sigma \right) =0\right\} $
and $\mathrm{D}\left( x,\sigma ,t\right) =\left( x,t\mathrm{\mathrm{d}}%
_{x}\Phi \left( x,\sigma \right) \right) .$ This implies\ that for any $x\in
X,\ $set $Z\left( x\right) =\left\{ \sigma ;\Phi \left( x,\sigma \right)
=0\right\} $ is diffeomorphic to the sphere $\mathrm{S}^{n-1}.$ Points $%
x,y\in X$ are called \textit{conjugate\ }for a generating function $\Phi $%
\textit{,} if $x\neq y,\ \Phi \left( x,\sigma \right) =\Phi \left( y,\sigma
\right) =0$ and \textrm{d}$_{\sigma }\Phi \left( x,\sigma \right) \parallel 
\mathrm{d}_{\sigma }\Phi \left( y,\sigma \right) $ for some $\sigma \in
\Sigma .$ Under condition (\textbf{I)} and condition (\textbf{II}): there
are no conjugate points, the integral%
\[
Q_{n}\left( x,y\right) =\int_{Z\left( y\right) }\left( \Phi \left( x,\sigma
\right) -i0\right) ^{-n}\frac{\mathrm{d}\Sigma }{\mathrm{d}_{\sigma }\Phi
\left( y,\sigma \right) }
\]%
is well defined for any $x,y\in X,$ $y\neq x$.

\begin{theorem}
\label{B}Let $\mathrm{d}X$ be the volume form of a Riemannian metric $%
\mathrm{g}$ on $X$ and $\Phi $ be a generating function satisfying (\textbf{I%
}),(\textbf{II})\textrm{\ }and condition (\textbf{III}): 
\begin{equation}
\func{Re}i^{n}Q_{n}\left( x,y\right) =0\ \text{for all }x,y\in X\text{ such
that}\ x\neq y.  \label{6}
\end{equation}%
For any odd $n,$ an arbitrary function $f\in \mathcal{C}^{n-1}\left(
X\right) $ with compact support can be reconstructed from data of the
Funk-Radon transform by%
\begin{equation}
f\left( x\right) =\frac{1}{2\mathrm{j}^{n-1}D_{n}\left( x\right) }%
\int_{\Sigma }\delta ^{\left( n-1\right) }\left( \Phi \left( x,\sigma
\right) \right) \mathrm{M}_{\Phi }f\left( \sigma \right) \mathrm{d}\Sigma .
\label{9}
\end{equation}%
For even $n,$ any function $f\in \mathcal{C}^{n-1+\varepsilon }\left(
X\right) ,$ is recovered by%
\begin{equation}
f\left( x\right) =\frac{\left( n-1\right) !}{\mathrm{j}^{n}D_{n}\left(
x\right) }\int_{\Sigma }\frac{\mathrm{M}_{\Phi }f\left( \sigma \right) 
\mathrm{d}\Sigma }{\Phi \left( x,\sigma \right) ^{n}},  \label{8}
\end{equation}%
where for any $n$ 
\[
D_{n}\left( x\right) =\frac{1}{\left\vert \mathrm{S}^{n-1}\right\vert }%
\int_{Z\left( x\right) }\frac{1}{\left\vert \nabla _{x}\Phi \left( x,\sigma
\right) \right\vert _{\mathrm{g}}^{n}}\frac{\mathrm{d}\Sigma }{\mathrm{d}%
_{\sigma }\Phi \left( x,\sigma \right) }.
\]%
The integrals (\ref{8}) and (\ref{9}) converge to $f$ uniformly on any
compact set $K\subset X.$
\end{theorem}

See \cite{Pa16} for a proof. The singular integrals like (\ref{8}) and (\ref%
{9}) are defined as follows%
\begin{eqnarray*}
\int \frac{\mathrm{\omega }}{\Phi ^{n}} &=&\frac{1}{2}\left( \int \frac{%
\mathrm{\omega }}{\left( \Phi -i0\right) ^{n}}+\int \frac{\mathrm{\omega }}{%
\left( \Phi +i0\right) ^{n}}\right) , \\
\int_{\Sigma }\delta ^{\left( n-1\right) }\left( \Phi \right) \mathrm{\omega 
} &=&\left( -1\right) ^{n-1}\frac{\left( n-1\right) !}{2\pi i}\left( \int 
\frac{\mathrm{\omega }}{\left( \Phi -i0\right) ^{n}}-\int \frac{\mathrm{%
\omega }}{\left( \Phi +i0\right) ^{n}}\right)
\end{eqnarray*}%
for any smooth $n$-form $\mathrm{\omega .}$

\section{Reconstructions}

\begin{theorem}
Let $E^{n+1}$ be an affine space with an invariant volume form $\mathrm{d}V$
and $X$ be a smooth hypersurface in $E^{n+1}$ (occasionally not closed) with
a Riemannian metric $\mathrm{g}$. Let $\Sigma $ be an ellipsoid in $E^{n+1}$
such that \newline
(\textbf{E}): any line that meets $X$ at least twice\ or is tangent to $X$
does not touch $\Sigma .$ Then for any odd $n,$ any function $f\in \mathcal{C%
}_{0}^{n-1}\left( X\right) $ can be recovered from integrals 
\[
\mathrm{M}_{\Phi }f\left( \sigma \right) =\int \delta \left( \left\langle
x-\sigma ,\nabla q\left( \sigma \right) \right\rangle \right) f\left(
x\right) \mathrm{d}_{\mathrm{g}}X=\int_{Z\left( \sigma \right) }f\left(
x\right) \frac{\mathrm{d}_{\mathrm{g}}X}{\left\langle \mathrm{d}x,\nabla
q\left( \sigma \right) \right\rangle },\ \sigma \in \Sigma 
\]%
by%
\[
f\left( x\right) =\frac{1}{2\mathrm{j}^{n-1}D_{n}\left( x\right) }%
\int_{\Sigma }\delta ^{\left( n-1\right) }\left( \left\langle x-\sigma
,\nabla q\left( \sigma \right) \right\rangle \right) \mathrm{M}_{\Phi
}f\left( \sigma \right) \mathrm{d}\Sigma .
\]%
For any even $n,$ an arbitrary function $f\in \mathcal{C}_{0}^{n-1+%
\varepsilon }\left( X\right) $ can be reconstructed by%
\[
f\left( x\right) =\frac{\left( n-1\right) !}{\mathrm{j}^{n}\Delta _{n}\left(
x\right) }\int_{\Sigma }\frac{\mathrm{M}_{\Phi }f\left( \sigma \right) }{%
\left\langle x-\sigma ,\nabla q\left( \sigma \right) \right\rangle ^{n}}%
\mathrm{d}\Sigma .
\]%
For any $n,$ $Z\left( \sigma \right) =\left\{ x;\left\langle x-\sigma
,\nabla q\left( \sigma \right) \right\rangle =0\right\} ,\ \mathrm{d}\Sigma =%
\mathrm{d}V/\mathrm{d}q$ and%
\[
\Delta _{n}\left( x\right) \doteqdot \frac{1}{\left\vert \mathrm{S}%
^{n-1}\right\vert }\int_{Z\left( x\right) }\frac{1}{\left\vert \nabla
_{x}\Phi \right\vert _{\mathrm{g}}^{n}}\frac{\mathrm{d}\Sigma }{\mathrm{d}%
_{\sigma }\left\langle x-\sigma ,\nabla q\left( \sigma \right) \right\rangle 
},\ x\in X.
\]%
The integrals converges uniformly on any compact set in $X.$
\end{theorem}

\textbf{Remark 1.} Here $\mathrm{d}_{\mathrm{g}}X$ is the volume form of the
Riemannian metric $\mathrm{g\ }$on $X$ and $\left\vert \cdot \right\vert _{%
\mathrm{g}}$ is the Riemannian norm of a covector.

\textbf{Remark 2. }Generating function $\Phi $ can be replaced by 
\begin{equation}
\Phi ^{\prime }\left( x,\sigma \right) =\left\langle x-e,\nabla q\left(
\sigma \right) \right\rangle -r,\ r=2-2q\left( e\right)   \label{10}
\end{equation}%
where $e$ is the center of $\Sigma .$ It\ coincides with $\Phi $ on $X\times
\Sigma $ and is linear in $\sigma .$ The ellipsoid can be replaced by an
arbitrary hyperboloid $H$ in $E^{n+1}$ and the volume form $T^{\ast }\left( 
\mathrm{d}V\right) $ where $T$ is a projective transform of the ambient
projective space $P^{n+1}$ such that $T\left( H\right) $ is an ellipsoid.

\textbf{Remark 3. }If $Y$ is a closed and convex manifold in $E^{n+1}$, and
the cam is inside of $Y,\ $then for any $\tau \in \Sigma ,$ the manifold $%
X=\left\{ x\in Y:\Phi \left( x,\tau \right) >0\right\} $ fulfils (\textbf{E)}%
. If $Y$ is a sphere with the center at a one-point cam then $X$ is a
hemisphere and Theorem 2 provides inversion of the Funk theorem \cite{Funk}.
If $X$ is a hyperplane and the cam is a point in $E^{n+1}$ or at infinity
Theorem 2 is equivalent to the Radon inversion theorem. In the latter case
to fulfil (\textbf{I) }one need to take each hyperplane through the cam
point two times with the opposite conormal vectors.

Theorem 2 was obtained in \cite{Pa16} for the case $X\ $is a subset of the
sphere \textrm{S}$^{n}$ and the cam is a sphere in the inside $\mathrm{S}^{n}
$. This result with one-point cam was considered also in \cite{Sal16}. Note
that for arbitrary $X$ and one-point cam $\left\{ e\right\} ,$ Theorem 2 is
reduced to Funk's result by the central projection of $X$ to the unit sphere 
$S^{n}$ with the center $e.$

\section{Proof}

The reconstruction formulas as above are invariant with respect to any
affine transformation. Therefore we can assume that $\Sigma $ is a sphere
(which makes some geometric arguments more obvious). The function $\Phi $
generates the family of hyperplane sections of $X$ since any hyperplane $H$
tangent to $\Sigma $ can be written in the form $H_{\sigma }=\left\{
x;\left\langle x-\sigma ,\nabla q\left( \sigma \right) \right\rangle
=0\right\} $ for some $\sigma \in \Sigma \mathrm{.}$ Now we check that $\Phi
\left( x,\sigma \right) =\left\langle x-\sigma ,\nabla q\left( \sigma
\right) \right\rangle $ satisfies conditions (\textbf{I}), (\textbf{II}) and
(\textbf{III}) as in Sect.2.

\newpage

\begin{lemma}
$\Phi $ satisfies (\textbf{I}).
\end{lemma}

\textit{Proof. }We have $\mathrm{d}_{x}\Phi \neq 0$ on $Z$ since of (\textbf{%
E). }For any point $x\in X$ and any covector $v\in T_{x}^{\ast }\left(
X\right) ,$ $v\neq 0,$ there exists one and only one hyperplane $H_{\sigma }$
such that $x\in H_{\sigma }$ and $v=t\mathrm{d}_{x}\Phi \left( x,\sigma
\right) $ on $T_{x}^{\ast }\left( X\right) $ for some $t>0$\textbf{. }It
follows that the map $\mathrm{D}_{X}$ is bijective. We prove that $\mathrm{D}%
_{X}$ is a local diffeomorphism. This condition can be written in the form%
\begin{equation}
\det J_{\xi ,\tau }\left( x,\sigma \right) \neq 0,\ \left( x,\sigma \right)
\in Z,  \label{2}
\end{equation}%
where%
\[
J_{\xi ,\tau }=\left( 
\begin{array}{cc}
^{t}\nabla _{\xi }\Phi  & \nabla _{\xi }\nabla _{\tau }\Phi  \\ 
0 & \nabla _{\tau }\Phi 
\end{array}%
\right) 
\]%
and $\xi ,\tau $ are arbitrary local systems of coordinates on $X$ and $%
\Sigma $ respectively. Let $T=\left( t,t_{0}\right) $ be a $n+1$-vector such
that $TJ=0$ where%
\[
J_{\xi ,\tau }=\left( 
\begin{array}{cc}
\left\langle \partial x/\partial \xi ,\nabla q\right\rangle  & \left\langle
\left( \partial x/\partial \xi \times \partial \sigma /\partial \tau \right)
,\nabla ^{2}q\left( \sigma \right) \right\rangle  \\ 
0 & \left\langle \left( x-\sigma \right) \times \partial \sigma /\partial
\tau ,\nabla ^{2}q\left( \sigma \right) \right\rangle 
\end{array}%
\right) 
\]%
and $\left\langle \partial \sigma /\partial \tau ,\nabla q\left( \sigma
\right) \right\rangle =0$ since $q$ is constant on $\Sigma .$ Equation $TJ=0$
is equivalent to%
\begin{eqnarray}
\left\langle \left\langle t,\partial x/\partial \xi \right\rangle ,\nabla
q\right\rangle  &=&0,\   \label{3} \\
\left\langle \left\langle t,\partial x/\partial \xi \right\rangle
+t_{0}\left( x-\sigma \right) \times \partial \sigma /\partial \tau
_{j},\nabla ^{2}q\left( \sigma \right) \right\rangle  &=&0,\ j=1,...,n.
\label{4}
\end{eqnarray}%
Vector $\left\langle t,\partial x/\partial \xi \right\rangle $ is tangent to 
$X$ and (\ref{3}) means that it is also tangent to $\Sigma \mathrm{\ }$at$\
\sigma .$ Vector $x-\sigma $ is also tangent to $\Sigma $ since of $\Phi
\left( x,\sigma \right) =0.$ Therefore there exist constants $c_{1},..,c_{n}$
such that%
\[
\theta \doteqdot c_{1}\partial \sigma /\partial \tau _{1}+...+c_{n}\partial
\sigma /\partial \tau _{n}=\left\langle t,\partial x/\partial \xi
\right\rangle +t_{0}\left( x-\sigma \right) .
\]%
Taking the corresponding linear combination of equations (\ref{4}) we get%
\[
\left\langle \theta \times \theta ,\nabla ^{2}q\left( \sigma \right)
\right\rangle =0
\]%
which implies $\theta =0$ since the form $\nabla ^{2}q$ is strictly
positive. It follows that $\left\langle t,\partial x/\partial \xi
\right\rangle +t_{0}\left( x-\sigma \right) =0$ which\ implies that both
terms vanish since the first one is tangent to $X$ and the second one is
transversal to $X.$ Finally $t=0$, $t_{0}=0$ and $T=0$ which completes the
proof of (\ref{2}) and of the Lemma. $\blacktriangleright $

Condition (\textbf{II}). Check that generating\ function $\Phi $ coincides
with (\ref{10}).  This follows from 
\begin{equation}
\Phi ^{\prime }\left( x,\sigma \right) -\Phi \left( x,\sigma \right)
=\left\langle \sigma -e,\nabla q\left( \sigma \right) \right\rangle
-r=2\left( q\left( \sigma \right) -q\left( e\right) \right) -r=0  \label{5}
\end{equation}%
since $q\left( \sigma \right) -q\left( e\right) $ is a quadratic form of $%
\sigma -e,$ $\sigma \in \Sigma $. Suppose that this condition violates for
some points $x,y\in X.$ We have then $a\left\langle x-e,\nabla ^{2}q\left(
\sigma \right) \right\rangle =b\left\langle y-e,\nabla ^{2}q\left( \sigma
\right) \right\rangle $ for some vector $\left( a,b\right) \neq \left(
0,0\right) $ and a point $\sigma \in \Sigma .$ This implies that $a\left(
x-e\right) =b\left( y-e\right) $ since the matrix $\nabla ^{2}q$ is
nonsingular. This yields that $x,y,$ and $e$ belong to one line. This line
crosses the cam which is impossible since of (\textbf{E)}.

\begin{lemma}
Function $\Phi $ fulfils (\textbf{III})\textbf{.}
\end{lemma}

\textit{Proof. }We are going to show that integral%
\begin{equation}
Q_{n}\left( x,y\right) =\func{Re}i^{n}\int_{Z\left( y\right) }\left( \Phi
\left( x,\sigma \right) -i0\right) ^{-n}\frac{\mathrm{d}\Sigma }{\mathrm{d}%
_{\sigma }\left\langle y-\sigma ,\nabla q\left( \sigma \right) \right\rangle 
}  \label{1}
\end{equation}%
vanishes for all $x,y\in X,\ y\neq x.$ We have%
\[
\Phi \left( x,\sigma \right) =\Phi \left( x,\sigma \right) -\Phi \left(
y,\sigma \right) =\left\langle x-y,\nabla q\left( \sigma \right)
\right\rangle 
\]%
on $Z\left( y\right) .$ The right hand side does not change its sign if and
only if the point $x$ is contained in the convex closed cone bounded by the\
lines through points $y$ that are tangent to $Z\left( y\right) .$ It is not
the case since of (\textbf{E}). Therefore $\Phi \left( x,\sigma \right) $
does change its sign on $Z\left( y\right) .$ By (\ref{5})%
\[
\mathrm{d}_{\sigma }\left\langle y-\sigma ,\nabla q\left( \sigma \right)
\right\rangle -\left\langle y-e,\mathrm{d}_{\sigma }\nabla q\left( \sigma
\right) \right\rangle =-\mathrm{d}\left( \sigma -e,\nabla q\left( \sigma
\right) \right) =-\mathrm{d}\left( q\left( \sigma \right) -q\left( e\right)
\right) =0
\]%
on $\Sigma $ since $q\left( \sigma \right) =1.$ The volume form in (\ref{1})
equals 
\[
\frac{\mathrm{d}\Sigma }{\mathrm{d}_{\sigma }\left\langle y-\sigma ,\nabla
q\left( \sigma \right) \right\rangle }=\frac{\mathrm{d}V}{\mathrm{d}q\wedge
\left\langle y-e,\mathrm{d}_{\sigma }\nabla q\left( \sigma \right)
\right\rangle }.
\]%
Choose affine\ coordinates $\sigma =A\xi +e$ on $E^{n+1}$ where $A$ is the
diagonal matrix such that $2q\left( A\xi +e\right) =\left\vert \xi
\right\vert ^{2}.$ Then $\mathrm{d}V=\mathrm{\det }A\mathrm{\ d}\xi
_{1}\wedge ...\wedge \mathrm{d}\xi _{n+1},\ \mathrm{d}q=\sum \xi _{i}\mathrm{%
d}\xi _{i}$ and $\left\langle y-e,\mathrm{d}\nabla q\left( \sigma \right)
\right\rangle =\left\langle s,\mathrm{d}\xi \right\rangle $ for some vector $%
s\in E^{n+1}.$ This yields 
\[
\frac{\mathrm{d}V}{\xi \mathrm{d}\xi \wedge \left\langle s,\mathrm{d}\xi
\right\rangle }=\frac{\Omega _{n}}{\left\langle s,\mathrm{d}\xi
\right\rangle }=\left\vert s\right\vert ^{-1}\Omega _{n-1}
\]%
up to the factor $\mathrm{\det }A.$ Here $\Omega _{k}$ denotes the volume
form of the euclidean $k$-sphere $\mathrm{S}^{k}$. Finally, we apply \cite%
{Pa16}Theorem A.20 to $\Phi $ and to the sphere $Z\left( y\right) \cong 
\mathrm{S}^{n-1}$ which implies vanishing of  $Q_{n}\left( x,y\right) .$ $%
\blacktriangleright $

Application of Theorem \ref{B} completes the proof of Theorem 2 for any
nondegenerated ellipsoid. In the case of one-point cam$\ \left\{ e\right\} \ 
$one can take the generating function $\Phi \left( x,\sigma \right)
=\left\langle x-e,\sigma \right\rangle ,$ $\sigma \in \mathrm{S}^{n}$ and
follow the above arguments. $\blacktriangleright $

\end{document}